\documentclass[12pt]{article}
\usepackage{amsmath, amssymb, amsthm}
\usepackage{enumerate}
\newtheorem{theorem}{Theorem}[section]
\newtheorem{lemma}{Lemma}[section]
\newtheorem{proposition}{Proposition}[section]

\newtheorem*{theorema}{$T(1)$ Theorem}
\newcommand{\taf}{{\hskip 5pt} $\blacksquare$
                  \renewcommand{\qedsymbol}{}}
\begin{document}
\title{Characterizations of the Hardy Space $H^1$\\ and BMO}
\author{Wael Abu-Shammala, Ji-Liang Shiu,\\ and Alberto Torchinsky}
\date{}
 \maketitle
\begin{abstract}
We describe the spaces $H^1(R)$ and BMO$(R)$ in terms of their
closely related, simpler  dyadic and two-sided counterparts. As a
result of these characterizations we establish when a bounded
linear operator defined on dyadic or two-sided $H^1(R)$ into a
Banach space has a continuous extension to $H^1(R)$ and when a
bounded linear operator that maps a Banach space  into dyadic or
two-sided  BMO$(R)$ actually maps continuously into BMO$(R)$.
\end{abstract}
\section{Introduction}
In this paper  we seek to elucidate the role  simple atoms, such
as the Haar system, play in the theory of the Hardy space
$H^1(R)=H^1$. It becomes quickly apparent that the Haar system, or
more generally dyadic atoms, do not suffice to span $H^1$.
However, the fact  that arbitrary atoms can be written as the sum
of at most three atoms, two dyadic and a special atom, makes it
possible to gain a greater insight into the structure of $H^1$ and
its dual, BMO. We pass now to describe the specific results.

 By the Hardy space $H^1(R)=H^1$ we mean the
collection of those integrable functions  which admit an atomic
decomposition in terms of $L^{\infty}$ atoms.  Recall that a
compactly supported function $a$ with vanishing integral is an
$L^{\infty}$ atom, or plainly an atom, with defining interval $I$,
if
\[{\text{supp$(a)\subseteq I$}},\quad |a(x)|\le \frac{1}{|I|}\,,
\quad {\text{and }}\, \int_I a(x)\,dx=0\,.\]
$H^1$ is then the
Banach space consisting of those $f$'s such that
\[H^{1}=\Big\{f=\sum_1^{\infty}\lambda_j\,a_j \, :
\sum_1^{\infty}|\lambda_j|<\infty\,\Big\},\]  where the
convergence is in the sense of distributions as well as in $L^1$,
the $a_j$'s are $L^{\infty}$ atoms, and the atomic norm is given
by
\[\|f\|_{H^1}=\inf\Big\{\sum_1^{\infty}|\lambda_j|:
f=\sum_1^{\infty}\lambda_j\,a_j \,\Big\}.\]
 For these, and all other well-known basic facts
 used throughout the article, see \cite {GCRdF},
 \cite {stein} and \cite {torchinsky}.

$H^1_d(R)=H^1_d$, or dyadic $H^1$, is obtained by restricting the
defining intervals to be dyadic. Clearly
$\|f\|_{H^1}\le\,\|f\|_{H^1_d}$, and the inclusion
   $H^1_d\subset H^1$ is strict.
 As for $\overline{H^1_d}$, the closure of $H^1_d$ in
$H^1$, it turns out to be the space
\[H^A=\Big\{f\in H^1:\int_0^{\infty}f(x)\,dx=0\Big\}.\]
More to the point: if $\{H_I\}$ denotes the $H^1$ normalized Haar
system indexed by the dyadic intervals $I$ of $R$, then
$\overline{{\text {sp}}}\{H_I\}$, the closed span of the Haar
system in $H^1$, is also $H^A$.

$H^A$  is not a convenient space -- for instance, $f\in H^A$ does
not imply that $f\chi_{[0,\infty)}$ belongs to $H^1$ -- and we are
led to introduce two-sided $H^1$, or $H^1_{2s}(R)=H^1_{2s}$. This
is the space of $H^1$ functions with atomic decompositions in
terms of atoms whose defining interval lies on either side of the
origin. Endowed with the atomic norm, $H^1_{2s}$ is a proper
subspace of $H^A$ continuously included in $H^1$. Also,  if $f\in
H^1_d$, $\|f\|_{H^1_{2s}}\,\le\,\|f\|_{H^1_d}$, and $\overline
{\text{sp}}^{\raisebox{-.925pt}{\hspace*{1pt}$\scriptstyle
2s$}}\{H_I\}$, the closed span of the Haar system in $H^1_{2s}$,
is $H^1_{2s}$,  see \cite{shiu}.

As for the atoms themselves, we note that each   atom can be
written as the sum of at most three atoms, two dyadic and a
special atom, see \cite{fridli}. This allows us to identify $H^1$
as the sum of $H^1_d$ and a space generated by special atoms and
to establish the boundedness of linear operators from $H^1$ taking
values in a Banach space $X$ given that they map $H^1_d$ into $X$.
Similar results hold for $H^A$ and $H^1_{2s}$.

The space of functions of bounded mean oscillation, BMO$(R)$=BMO,
consists of those locally integrable functions $\varphi$  with
\[\|\varphi\|_*=\sup_I \frac{1}{|I|}\int_I \,|\,\varphi(x)-
\varphi_I\,|\, dx<\infty\,,\] where
$\varphi_I=\frac1{|I|}\int_I\,\varphi (y)\,dy$ is the average of
$\varphi$ over $I$, and the $\sup$ is taken over all bounded
intervals $I$. Modulo constants, $({\rm BMO}, \|\cdot\|_*)$ is a
Banach space.

${\rm BMO}_d(R)={\rm BMO}_d$, or dyadic BMO, is defined by
restricting the intervals above  to be dyadic. The $\sup$ is now
denoted by $\|\cdot\,\|_{*,d}$, and modulo two constants (one for
each side of the origin), $({\rm BMO}_d, \|\cdot\,\|_{*,d})$ also
becomes a Banach space.

Finally, when the above $\sup$ is restricted to those  $I$'s that
lie on either side of the origin,  the sup is denoted
$\|\,\cdot\,\|_{2s}$, and we have two-sided BMO or ${\rm
BMO}_{2s}(R)={\rm BMO}_{2s}$. Modulo two constants, $({\rm
BMO}_{2s},\|\,\cdot\,\|_{2s})$ becomes a Banach space.

BMO is the dual of $H^1$, ${\rm BMO}_d$ is the dual of $H^1_d$,
 and below note that ${\rm BMO}_{2s}$ is the dual of $H^1_{2s}$ and
identify the dual of $H^A$. Each of the decompositions of $H^1$
suggests a characterization of BMO in terms of these dual spaces.
These new characterizations of BMO in turn allow us to work in
dyadic and dyadic-like settings, and
 provide us with an effective   way to pass from ${\rm BMO}_d$,
 and ${\rm BMO}_{2s}$,  to
BMO.  We also establish when a bounded linear operator that maps a
Banach space $X$ into ${\rm BMO}_d$ actually maps $X$ into BMO.

 The paper is
organized as follows. Section two is devoted to  $H^1$,  in
section three we discuss BMO, and we  conclude the paper in
section four by  pointing out how these results can be extended to
$R^n$, $n>1$.
\section{The Hardy Space $H^1(R)$}
\subsection*{The Haar System}
A dyadic interval $I$ is an interval of the special form
\[I=I_{n,k}=[\, k 2^n, (k+1)2^n),\]
where $k$ and $n$ are arbitrary integers, positive, negative or
$0$. Note that $I=I_L \cup I_R$, where the left half $I_L$ and the
right half $I_R$ of $I$ are also dyadic.

For each dyadic interval $I$,  the $H^1$ normalized Haar function
$H_I$ is given by
\[ H_I(x)= \frac1{{|I|}} \chi_{I_L}(x) -
 \frac{1}{{|I|}} \chi_{I_R}(x)\,.\]

Finally, throughout the paper $b$ will denote the function
\[b(x)=\tfrac12\left[\chi_{(-1,0)}(x)-\chi_{(0,1)}(x)\right].\]

To describe how $H^1_d$ fits in $H^1$, let $a$ be an atom.
 If the origin is not an interior point
of the defining interval of $a$, $a$ is a multiple of a dyadic
atom. On the other hand, if the origin is interior to the defining
interval of $a$, let
\[a(x) = \Big[a(x)+2\Big(\int_0^{\infty}a(y)\,dy\Big)\,b(x)\Big] -
2\Big(\int_0^{\infty}a(y)\,dy\Big)\,b(x)\,.\] Since $a$ has
vanishing integral the first function above is a linear
combination of two dyadic atoms with defining intervals on
opposite sides of the origin, and the second is a multiple of the
fixed function $b$. Thus $H^1_d$ is of codimension one in $H^1$,
and since $\overline{H^1_d}\subset H^A$, actually
$\overline{H^1_d}=H^A$.

To show that the same is true for the closed span of the Haar
system in $H^1$, we make use of the following observation. Its
proof is  left to the reader.
\begin{lemma} Suppose a locally integrable function
$\varphi$ satisfies
\[ \int_R H_{I}(x) \varphi(x)\, dx = 0\quad {\text{ for all
dyadic intervals $I$}}\,.\] Then for some constants $c,d$,
\[
\varphi(x)=d\,\chi_{(-\infty, 0)}(x) + c \,\chi_{(0,
\infty)}(x)\,.\]
\end{lemma}
We  then have,
\begin{theorem}
The closed span of the Haar system in $H^1$ is $H^A$.
\end{theorem}
\begin{proof}
It suffices to show that $\overline
{\text{sp}}^{\raisebox{-.925pt}{\hspace*{1pt}$\scriptstyle
d$}}\{H_I\}$, the closed span of the Haar system in $H^1_d$, is
$H^1_d$. Let $L$ be a bounded linear functional  on $H^1_d$ that
vanishes on the $H_I$'s.  Then there is  $\varphi \in {\rm
BMO}_d(R)$ with $\|\varphi \|_{*} \sim \|L\|$ such that for
compactly supported $f\in H^1_d$, $L(f) = \int_R f(x) \varphi(x)\
dx$.  Now, since $L$ vanishes on the $H_{I}$'s,  by Lemma 2.1,
$\varphi(x)=d\,\chi_{(-\infty, 0)}(x) + c \,\chi_{(0,
\infty)}(x)$, and  consequently for those $f$'s
\[L(f)= d \int_{-\infty}^0 f(x) \, dx + c \int_0^{\infty}
f(x)\, dx=0\,.\] Thus  $L$ is the zero functional and we have
finished.\taf
\end{proof}
The reader will have no difficulty in establishing the
quantitative version of Theorem 2.1 in terms of $d(f,H^A)$, the
distance of $f$ in $H^1$ to $H^A$.
\begin{proposition}Let $f\in H^1$. Then
$d(f,H^A)\sim \Big|\int_0^{\infty}f(x)\,dx\Big|$.
\end{proposition}
\subsection*{Structure of Atoms}
Since dyadic atoms do not suffice to represent  $H^1$ functions,
we consider how far apart dyadic atoms are from arbitrary atoms.
The answer is that
 an arbitrary $H^1$ atom can be expressed as a sum of
at most three atoms, two dyadic and a special atom, see
\cite{fridli}.
\begin{lemma}
 Let $a$ be an $H^1$ atom. Then there are at most
three atoms $a_L, a_R, b_{n,k}$, such that
\begin{enumerate}[\upshape 1.]
\item  $a_L$ and $a_R$ are dyadic atoms.
\item For some integers
$n,k$,
\[b_{n,k}(x)= \frac {1}{2^{n+1}}\left[ \chi_{[(k-1)2^n,\,
k\, 2^n]}(x)-  \chi_{[k\, 2^n, (k+1) 2^n]}(x)\right].\] \item
$a=c_1a_L+c_2a_R+c_3b_{n,k}$, where $|c_1|,|c_2|, |c_3|\le 4\,$.
\end{enumerate}
\end{lemma}
\begin{proof}
  Let $I$ be the defining interval for $a$, and let $n$
be the integer such that $2^{n-1}\le |I| <2^n$ and $k$ the integer
such that  $I\subset [\,k\,2^n, (k+1)2^n]$. Set now $a_L$ equal to
\[
a_L(x) =
 \begin{cases} \displaystyle  \frac{1}{4}\Big( a(x)
-\frac{1}{2^n}\int_{[(k-1)2^n,\,k\,2^n)}
 a(y)\, dy\Big),
& x\in [(k-1)2^n,k\,2^n),\\[5pt]
\displaystyle 0, & \mbox{otherwise.}
\end{cases}
\]
Since $a$ is an atom with defining interval $I$ it readily follows
that
\[ \|a_L\|_{\infty} \le \frac{1}{4}\Big(\frac1{|I|} +
\frac1{|I|}\Big)\le \frac1{2\,|I|}\le 2^{-n}\,.
\]
Furthermore, since $a_L$ is supported in $[(k-1)2^n,\,k\,2^n]$ and
has integral $0$,  $a_L$ is a dyadic atom.

Similarly, set $a_R$ equal to
\[
a_R(x) =\begin{cases} \displaystyle \frac{1}{4}\Big( a(x)
-\frac{1}{2^n}\int_{[k\,2^n,\,(k+1)2^n)}a(y)\, dy\Big),
& x\in [k\,2^n,(k+1)2^n),\\[5pt]
\displaystyle 0, & \mbox{otherwise.}
\end{cases}
\]
$a_R$ is supported in $[k\,2^n,\,(k+1)2^n]$, $\|a_R\|_{\infty} \le
2^{-n}$, and has integral $0$,   so $a_R$ is also a dyadic atom.

Finally put
\[b_{n,k}(x)= \frac {1}{2^{n+1}}\Big[ \chi_{[(k-1)2^n,\,
k\, 2^n]}(x)-  \chi_{[k\, 2^n,\, (k+1) 2^n]}(x)\Big].\] Since $a$
has vanishing integral and $I\subset[(k-1)2^n, (k+1)2^n]$  it is
also true that $-\int_{[(k-1)2^n,\, k\,2^n]} a(y)\, dy =
\int_{[k\,2^n,\, (k+1)2^n]} a(y)\, dy\,,$ and consequently, since
 $a(x)-4a_L(x)-4a_R(x)$ is equal to
\[
 \begin{cases}
\displaystyle \frac{1}{2^n}\int_{[(k-1)2^n,\,k\,2^n]}a(y)\, dy\,,
& x\in [(k-1)2^n,k\,2^n],\\[12pt]
\displaystyle \frac{1}{2^n}\int_{[k\,2^n,\,(k+1)2^n]}a(y)\, dy\,,
 & x\in [k\,2^n,(k+1)2^n],
\end{cases}
\]
we have
\[a(x)=4a_L(x)+4a_R(x)+2\Big(\int_{[(k-1)2^n,\,k\,2^n]}a(y)\ dy\Big)
b_{n,k}(x)\,.\] Thus the conclusion follows in this case with
$c_1=c_2=4$ and
\[|c_3|\le  2\, \int_{[(k-1)2^n,\,k\,2^n]}|a(y)|\ dy \le 2\,
\frac {2^n}{|I|}\le 4\,.{\hskip 10pt}\blacksquare\]
\renewcommand{\qedsymbol}{}
\end{proof}
\subsection*{$H^1$ as the sum of Banach Spaces}
As a first application of the decomposition of atoms we will show
that $H^1$ can be written as the sum of various Banach spaces. We
have already seen that $H^1=H^A+{\rm sp}\{b\}$ as linear spaces.
In fact, the reader will have no difficulty in verifying that
actually $H^1=H^A +\,{\text{sp}}\{h\}$, where $h$ is an arbitrary
function in $H^1\setminus H^A$, and $\|f\|_{H^1}\sim
\|f\|_{H^A+\,{\text{sp}}\{h\}}$.

Furthermore, for integers $n,k$, let $b_{n,k}$ denote the dyadic
dilations and integer translations of $b$, i.e.,   the collection
of   atoms  given by
 \[b_{n,k}(x)=
\tfrac{1}{2^{n+1}}\left[\chi_{[(k-1)2^n, k\,2^n]}(x)-
\chi_{[k\,2^n,(k+1)2^n]}(x)\right]\,.
\]
Note that the special atoms $b_{n,k}$ are multiples of dyadic
atoms if $k$ is odd, but not if $k$ is even. Also, if $k\ne 0$,
the support of $b_{n,k}$ lies on one side of the origin.

 Let $H^1_{\delta}(R)=H^1_{\delta}$ be the  space which consists of
the linear combinations
\[H^1_{\delta}=\Big\{\varphi=\sum_1^{\infty}\lambda_j\,{b_{n_j,k_j}} \,
: \sum_1^{\infty}|\lambda_j|<\infty\,\Big\}.\] Endowed with the
atomic norm
\[\|\varphi\|_{H^1_{\delta}}=\inf\Big\{\sum_1^{\infty}|\lambda_j|:
\varphi=\sum_1^{\infty}\lambda_j\,b_{n_j,k_j} \,\Big\},\]
$(H^1_{\delta},\|\cdot\|_{H^1_{{\delta}}})$ is a Banach space.
Observe that if $f\in H^1_{\delta}$, then $f\in H^1$ and
$\|f\|_{H^1}\le \|f\|_{H^1_{\delta}}$.

Similarly, when $k=0$ we denote the resulting space
$H^1_{\delta,0}$. It is clear that if $f\in H^1_{\delta,0}$, then
$f\in H^1_{\delta}$ and $\|f\|_{H^1}\le \|f\|_{H^1_{\delta}}\le
\|f\|_{H^1_{\delta,0}}$. $H^1_{\delta}$ and $H^1_{\delta,0}$ are the
spaces of special atoms alluded to above, see \cite{desouza}.

From Lemma 2.2 it readily follows that
\begin{proposition} $H^1 = H^1_d +H^1_{\delta}$, and
$\|f\|_{H^1}\sim \|f\|_{H^1_d +H^1_{\delta}}\,.$
\end{proposition}
The meaning of this decomposition is the following. The Haar
system, or more generally the dyadic atoms, divide the line in two
regions, $(-\infty,0\,]$ and $[\,0,\infty)$. To allow for the
information carried by a dyadic interval to be transmitted to an
adjacent dyadic interval, they must be  connected. The $b_{n,0}$'s
channel information across the origin and the remaining
$b_{n,k}$'s connect adjacent dyadic intervals that are not
subintervals of the same dyadic interval.

We also have
\begin{proposition} $H^1 = H^1_{2s}
+H^1_{\delta,0}$, and $\|f\|_{H^1}\sim \|f\|_{H^1_{2s}
+H^1_{\delta,0}}\,.$
\end{proposition}
 This characterization allows us
to identify the range of the projection mapping $P$ of $H^1$
functions $f$ into $Pf=\chi_{[0,\infty)}\,f$. It is  isomorphic to
$H^1_o$, the odd functions of $H^1$, and consists
 of those functions in $L^1(R^+)$  whose
Telyakovski\u{\i}\ transform also belongs to $L^1(R^+)$, see
\cite{fridli1}.

The reader will have no difficulty of verifying the following
observation, which is useful when considering mappings $T:H^1\to
X$.
\begin{proposition} Let $B=B_0 + B_1$, where $B_0,B_1$ are
Banach spaces, and assume  $T$ is a linear operator that maps
$B_0$ continuously into a Banach space $X$ with norm $\|T\|_0$.
 Then $T:B\to X$ with norm $\|T\|$
if and only if  $T:B_1\to X$ is bounded with norm $\|T\|_1$ and
\[\|T\| = \max\,(\,\|T\|_0\,,\|T\|_1\,).\]
\end{proposition}
We now apply Proposition 2.4 to Hardy type operators in the
setting $H^1=H^1_{2s} +H^1_{\delta,0}$. For $0\le\varepsilon\le 1$
let $\tau_{\varepsilon}$ be given by
\[\tau_{\varepsilon} f(x)=\frac1{|x|^{1-\varepsilon}}
\int_{-x}^xf(y)\,dy\,,\quad x\ne 0\,.
\]
We set $1/p=1-\varepsilon$, $1\le p\le \infty\,$, and consider
when $\tau_{\varepsilon}f$ is in $X=L^p(R)$. Since
$\tau_{\varepsilon}b_{n,0}=0$ for $b_{n,0}\in H^1_{\delta,0}$,
 the continuity on $H^1$ is
equivalent to that on $H^1_{2s}$. The case $\varepsilon=1$ is
trivial and merely states $\|\tau_{1}f\|_{\infty}\le \|f\|_{H^1}$.

For the remaining cases we begin by observing that for a two-sided
atom $a$ with defining interval $I$, $\tau_{\varepsilon} a$ is
also supported in $I$, and $\left( \int_{R}|\tau_{\varepsilon}
a(x)|^p\,dx \right)^{1/p}$ $\le \, \|\ln|\cdot|\,\|_{*,2s}^{1/p}$.
Let now $f\in H^1_{2s}$ have the atomic decomposition $f(x)=\sum_j
\lambda_j a_j\,$. Since the convergence also takes place  in
$L^1$, it readily follows that $\tau_{\varepsilon}
f(x)=\sum_j\lambda_j \,\tau_{\varepsilon} a_j (x)$, and thus by
Minkowski's inequality, upon taking the infimum over all possible
decompositions of $f$, we get
\[\Big(\int_{R}|\tau_{\varepsilon} f(x)|^p\,dx\Big)^{1/p}\le
\,\|\ln |\cdot|\,\|_{*,2s}^{1/p}\,\|f\|_{H^1_{2s}}.\] In short,
$\tau_{\varepsilon}$ maps $H^1$  into $L^p(R)$ with norm $\le
\|\ln |\cdot|\,\|_{*,2s}^{1/p}$. A similar reasoning applies to
the more general operators of Hardy type discussed in
\cite{GCRdF}, which include the Fourier transform.

Sublinear operators may be treated in a similar fashion. Consider,
for instance, $\widetilde{M}_{\varepsilon, d}$, the maximal
operator on $H^1_d$ given by
\[\widetilde{M}_{\varepsilon, d}f(x)=\sup_{x\in I}\frac1{|I|^{1-\varepsilon}}
\Big|\int_I f(y)\,dy\,\Big|,\] where $0\le\varepsilon<1$, and $I$
varies over the collection of dyadic intervals containing $x$, see
\cite{APC}. If $a$ is a dyadic atom with defining interval $I$,
\[\widetilde{M}_{\varepsilon, d}\,a(x)\le
\frac{1}{|I|^{1-\varepsilon}}\chi_I(x)\,,\] and consequently
\[\big\|\widetilde{M}_{\varepsilon, d}\,a\big\|_p\le\,c\,,\quad
0<\frac1p=1-\varepsilon\,.\] Thus $\widetilde{M}_{\varepsilon, d}$
is uniformly bounded on atoms and since it satisfies an
appropriate $\sigma$-sublinearity condition, it can be extended
continuously as a mapping from $H^1_d$ into $L^p$.  However, since
$\widetilde{M}_{\varepsilon, d}\,b(x)\sim \min(1,
|x|^{\varepsilon-1})$, $\widetilde{M}_{\varepsilon, d}$ only maps
the special  atoms into wk-$L^p$. So, $\widetilde{M}_{\varepsilon,
d}$ is bounded on $H^1_d$ but not on $H^1$.

On the other hand,  the  truncated version of this maximal
operator is better  behaved. For a positive integer $N$, let
\[ \widetilde{M}^N_{\varepsilon, d}f(x)=\sup_{x\in I}\frac1{|I|^{1-\varepsilon}}
\Big|\int_I f(y)\,dy\,\Big|\,,\] where $0\le\varepsilon<1$, and
$I$ varies over the collection of dyadic intervals  containing $x$
of size $2^{-N}\le |I|\le 2^N$. Then
$\widetilde{M}^N_{\varepsilon, d}$  also maps the $b_{n,k}$'s into
$L^p$, and this maximal function maps $H^1$ into $L^p$.
\subsection*{The duals of $H^1_{2s}$ and $H^A$}
An argument along by now familiar ideas allows us to identify the
dual of $H^1_{2s}$ as ${\rm BMO}_{2s}$. Indeed, we have
\begin{theorem} ${\rm BMO}_{2s}$ is the dual of $H^1_{2s}$.
More precisely,  for every $\varphi$ in ${\rm BMO}_{2s}$, the
functional
 $L_{\varphi}$ defined initially for bounded compactly supported
 functions $f\in H^1_{2s}$ by the integral $  L_{\varphi}(f)=
 \int_R f(x)\,\varphi(x)\,dx$
 has a bounded extension to $ H^1_{2s}$ with norm less than or
 equal to $c\,\|\varphi\|_{*,2s}$.

  Moreover, for any functional  $L\in (H^1_{2s})^*$, there is
   $\varphi\in {\rm BMO}_{2s}$, with norm $\|\varphi\|_{*,2s}\sim \|L\|$,
   such that $L(f)=L_{\varphi}(f)=\int_R f(x)\,\varphi(x)\,dx$
for every compactly supported  bounded $f\in H^1_{2s}$.
  \end{theorem}
We consider  the dual of $H^A$ next. Let $\varphi,\psi\in {\rm
BMO}$. We say that $\varphi\sim_2\psi$ if for some constants
$c,d$,
\[\varphi(x)-\psi(x)= d\chi_{(-\infty,0]}(x)
+c\chi_{[0,\infty)}(x)\, {\text{ a.e.}}\] Clearly $\sim_2$ is an
equivalent relation, and the norm for the element $\Phi\in B={\rm
{BMO}}/\sim_2$ with representative $\varphi\in {\rm BMO}$ is given
by
\[\|\Phi\|_B=\inf\{\|\psi\|_*: \varphi\sim_2 \psi\,\}\,.\]

Now,  since $(H^A)^*$ is isometrically isomorphic to
\[(H^1)^*/(H^A)^{\perp}={\rm{BMO}}/(H^A)^{\perp}\,,\]
we have  $(H^A)^* = B$.

This identification allows us to
 distinguish between  $H^1_{2s}$ and $H^A$. Indeed, since
$\chi_{(0,\infty)}(x)\,\ln x$ is in ${\rm BMO}_{2s}$ but not in
$B$, the inclusion $H^1_{2s}\subset H^A$ is strict.  We can also
exhibit  $f\in H^A\setminus H^1_{2s}$. For $n\ge 1$, let
$f_n(x)=2^n \chi_{[2^{-n},2^{-n+1}]}(x)$, and with
$L=\sum_{n=1}^{\infty} n^{-2}$, let $f(x)$ be the odd extension of
the function
\[L\chi_{[0,1]}(x)-\sum_{n=1}^{\infty}
\frac1{n^2}\,f_n(x)\,, \quad x>0\,.\] Then $f\in H^A$, and
$\|f\|_{H^1}\le 6\,L$. On the other hand, since
\[\int_0^{\infty} f_n(x)\,\ln x\,dx=-n\,\ln2+(2\ln2-1)\,,\]
it readily follows that $|\int_0^{\infty}f(x)\, \ln
x\,dx|=\infty$, and  $f\notin H^1_{2s}$. In other words, $f\in
H^A$ and its projection $Pf$ is an integrable function with
vanishing integral supported in $[0,\infty)$, but $Pf\notin H^1$.
However, if $g\in H^A$ vanishes for $x<0$, then $Pg\in H^1_{2s}$.
For, if $g$ has an atomic decomposition
$g(x)=\sum_j\lambda_j\,a_j(x)$, it also has the two-sided
decomposition
\[g(x)=\sum_j 4\,\lambda_j\,\frac14[a_j(x)+a_j(-x)]\chi_{[0,\infty)}(x)\,.\]
\section{Characterizations of BMO}
\subsection*{From $\rm{\bf BMO}_d$  to BMO} When restricted to
linear functionals, Proposition 2.4 suggests different
characterizations of BMO. We discuss the dyadic case first. Given
a BMO function $\varphi$,  consider the bounded linear functional
on $H^1$ induced by $\varphi$.  When acting on individual atoms,
two conditions, one for dyadic atoms and the other for special
atoms, must be satisfied for this functional to be bounded. The
 condition on the  dyadic atoms  suggests that
 $\varphi\in {\rm BMO}_d$, whereas the condition on the
 $b_{n,k}$'s, restates that the integral of $\varphi$ is in the Zygmund
class.

This motivates the following definition. For a locally integrable
function $\varphi$ let
 \[A(\varphi)=\sup_{n,k} \frac{1}{2^{n+1}} \Big|\int_{[(k-1)\,2^n,
 k\,2^n]}\varphi(x)\,dx- \int_{[k\,2^n,(k+1)\,2^n]}
 \varphi(x)\,dx\Big|\,,\]
 and put
 \[ \Lambda=\{\varphi\in {\rm BMO}_d: A(\varphi)<\infty\}\,, \quad
 \|\varphi\|_{\Lambda}=\max\Big(\|\varphi\|_{*,d}\,,
A(\varphi)\Big)\,.\]

Our next result describes how to pass from ${\rm BMO}_d$ to BMO.
\begin{theorem} ${\rm BMO}=\Lambda$. More precisely, if $\varphi\in {\rm BMO}$,
then $\varphi\in\Lambda$ and $\|\,\varphi\|_{\Lambda}\le
\|\,\varphi\|_{*}$. Also, if $\varphi\in\Lambda$, then $\varphi\in
{\rm BMO}$ and $\|\,\varphi\|_{*}\le c\,\|\,\varphi\|_{\Lambda}$.
\end{theorem}
\begin{proof}
It is clear that $\|\,\varphi\|_{\Lambda}\le c\,
\|\,\varphi\|_{*}$. Conversely, assume that $\varphi\in\Lambda$
and observe that for
 an atom $a=c_1a_L+c_2a_R+c_3b_{n,k}$, we have
\begin{align*}
\Big|\int_R a(&x)\varphi(x)\,dx\Big|
\\
 &\le 4\, \Big|\int_R
a_L(x)\varphi(x)\,dx\Big| + 4\, \Big|\int_R
a_R(x)\varphi(x)\,dx\Big|+4 \,\Big|\int_R
b_{n,k}(x)\varphi(x)\,dx\Big|\\[1pt]
 &\le 4\, \|\,\varphi\,\|_{*,d}+ 4\,
\|\,\varphi\,\|_{*,d}+4\,A(\varphi)\le\,12\,\|\,\varphi\|_{\Lambda}.
\end{align*}
 Suppose now that $f=\sum_1^{\infty}
\lambda_n a_n\in H^1$ is compactly supported and  bounded, and
that $\varphi\in\Lambda$ is bounded. Since $\lim_{N\to
\infty}\sum_1^N \lambda_n a_n=f$ in $L^1$,
  $\sum_1^N \lambda_n a_n \varphi$ converges
 to $f\varphi $ in $L^1$ as $N\to\infty$. Thus
$\lim_{N\to\infty}\sum_1^N\lambda_n\int_R a_n(x)\varphi(x)\,dx$
$=\int_R f(x)\varphi(x)\,dx$ and
\[\Big|\int_R f(x)\varphi(x)\,dx\Big|\le
\sum_1^{\infty}|\lambda_n| \Big|\int_R
a_n(x)\varphi(x)\,dx\Big|\le
c\,\sum_1^{\infty}|\lambda_n|\,\|\,\varphi\|_{\Lambda}.\] Since
the decomposition of $f$ is arbitrary, $\left|\int_R
f(x)\,\varphi(x)\,dx\right|\le c\,\|f\|_{H^1}\,
\|\,\varphi\|_{\Lambda}\,.$  To  show  that this estimate  also
holds for arbitrary $\varphi\in\Lambda$,  note that if $\varphi^k$
denotes the truncation of $\varphi$ at level $k$,
$\|\varphi^k\|_{\Lambda}\le c\,\|\varphi\|_{\Lambda}$ uniformly in
$k$. Thus
\[\Big|\int_R
f(x)\varphi^k(x)\,dx\Big|\le
c\,\|f\|_{H^1}\|\varphi^k\|_{\Lambda}\le
c\,\|f\|_{H^1}\|\varphi\|_{\Lambda}.\] Now, since $f$ has compact
support and $\varphi$ is locally integrable, $f\,\varphi$ is
integrable. Whence,   by the dominated convergence theorem,
$\int_R f(x)\varphi^k(x)\,dx$ tends to $\int_R
f(x)\varphi(x)\,dx$, and consequently,
\[\Big|\int_R f(x)\varphi(x)\,dx\Big|\le c\,
\|f\|_{H^1} \|\varphi\|_{\Lambda}\,.\] Finally, since
 \[ \|\varphi\|_{*}=\sup_{f\in H^1,\,\|f\|_{H^1}\le 1}
\Big|\int_R f(x)\varphi(x)\,dx\Big|,
 \]
 where $f$ is  compactly supported and bounded,
 $\|\varphi\|_{*}\le c\, \|\varphi\|_{\Lambda}\,.$
 \taf
 \end{proof}
As an application of the above characterization, the following
holds.
\begin{proposition}
Let  $T$ be a continuous linear operator  defined on a Banach
space $X$ which assumes values in ${\rm BMO}_d$ with norm
$\|T\|_d$. Then $T$ maps $X$ continuously into BMO if and only if
for all integers $n,k$,
\[\frac{1}{2^{n+1}} \Big|\int_{[(k-1)\,2^n,
 k\,2^n]}Tf(x)\,dx- \int_{[k\,2^n,(k+1)\,2^n]}
 Tf(x)\,dx\Big| \le M\,\|f\|_X\,,\]
and $\|T\| = \max\,(\,\|T\|_d\,,M)$.
 \end{proposition}
 \subsection*{Shifted BMO}
The process of averaging the translates of dyadic  BMO functions
leads to BMO,   and is an important tool in obtaining results in
BMO once they are known to be true in its dyadic counterpart,
${\rm BMO}_d$, see \cite{garnettjones}. It is also known that BMO
can be obtained as the intersection of ${\rm BMO}_d$ and one of
its shifted counterparts,  cf. \cite{tao_mei}.  These results
motivate the observations in this section.

Given a dyadic interval $I=[(k-1)2^{n},k\,2^n)$ of length $2^n$,
we call the interval
 $I'=[(k-1)2^n+2^{n-1}, k\,2^n+2^{n-1})$ the  shifted
 interval of  $I$ by its half-length. Clearly $|I'|=|I|$, and
 $I'=[(2k-1)2^{n-1},(2k+1)2^{n-1})$ is not dyadic.

Let $J=\{J_{n,k}\}$ be the collection of all dyadic shifted by
their half-length, $J_{n,k}=[(k-1)2^n, (k+1)2^n)$, all integers
$n,k$, and let ${\rm BMO}_{d^s}$ be the space consisting of those
locally integrable functions $\varphi$ such that
\[\|\varphi\|_{*,\,d^s}=\sup_{n,k}\frac{1}{|J_{n,k}|}\int_{J_{n,k}}
\left|\varphi(x)-\varphi_{J_{n,k}}\right|\,dx<\infty\,.\]

We then have \begin{theorem} ${\rm BMO}={\rm BMO}_d\cap {\rm
BMO}_{d^s}$.
\end{theorem}
\begin{proof} It is obvious that if $\varphi\in {\rm BMO}$, then
$ \|\varphi\|_{*,d}\,, \|\varphi\|_{*,d^s}\le \|\varphi\|_{*}$.

Conversely, it suffices to show that $\varphi\in \Lambda$. Since
$\varphi\in BMO_d$ it is enough to show that $A(\varphi)<\infty$.
For integers $n,k$, consider
 \begin{flushleft}
$\displaystyle {\frac{1}{2^{n+1}} \Big|\int_{[(k-1)\,2^n,
 k\,2^n]}\varphi(x)\,dx- \int_{[k\,2^n,(k+1)\,2^n]}
 \varphi(x)\,dx\Big|}$
 \end{flushleft}
\vspace{-6.5pt}
\begin{align*}&=\frac{1}{2^{n+1}} \Big|\int_{[(k-1)\,2^n,
 k\,2^n]}\Big[\varphi(x)-\varphi_{J_{n,k}}\Big]dx-
 \int_{[k\,2^n,(k+1)\,2^n]}
 \Big[\varphi(x)-\varphi_{J_{n,k}}\Big]dx\Big| \\[3pt]
&\le \frac{1}{2^{n+1}}\int_{[(k-1)2^n,(k+1)\,2^n]}
 \Big|\varphi(x)-\varphi_{J_{n,k}}\Big|\,dx\le
 \|\varphi\|_{*,d^s}\,,
 \end{align*}
 which implies that $A(\varphi)\le \|\varphi\|_{*,d^s}$, and
 we are done.
 \taf
 \end{proof}
 \subsection*{Further Characterizations of BMO}
We  further describe BMO in terms of the duals of the various
spaces describing $H^1$. From $H^1=H^A+{\rm sp}\{h\}$ it follows
that
 ${\rm BMO}=B\cap ({\text{sp}}\{h\})^*$, where $h\in H^1$ satisfies
 $\int_{[0,\infty)}h(y)\,dy\ne 0$. To fix ideas we pick $h=-b$ and
 introduce the equivalence relation $\sim_b$ in BMO as follows.
 We say that $\varphi \sim_b \psi$ if $\varphi-\psi=
 \eta$ for some  $\eta\in {\rm BMO}$ with
 $\int_R \eta(y)\,b(y)\,dy=0$. We endow these equivalence classes,
 which we denote by $B_b$,
 with the quotient norm,  and observe that the norm in BMO is
equivalent to the norm in $B\cap B_b$. It is possible, however,
to work with a simpler expression.
\begin{proposition} For a locally
integrable function $\varphi$, let
\[A_b(\varphi)=\Big|\int_R \varphi(y)\,b(y)\,dy\,\Big|.\] Let
$\varphi\in {\rm BMO}$ be the representative of $\Phi\in B$. Then
\[\|\varphi\|_*\sim
\max (\,\|\Phi\|_B, A_b(\varphi)\,).\]
\end{proposition}
\begin{proof}
If $\varphi\in {\rm BMO}$,  then clearly $\|\Phi\|_B\le
\,\|\varphi\|_*$. Also, $|A_b(\varphi)|\le \|\varphi\|_*
\|b\|_{H^1}$ $\le \|\varphi\|_*$.

As for the other inequality, we have $\|\varphi\|_*\le c\,\max
(\,\|\Phi\|_B, \|\tilde{\Phi}\|_{B_b})\,.$ Let $A=\int_R\varphi(y)
\,b(y)\,dy$, put $\psi(x)=2A\,b(x)$, and observe that $\psi\in
L^{\infty}$ and $\varphi(y)\sim_b\psi$. Then
$\|\tilde{\Phi}\|_{B_b}\le\,\|\psi\|_*\le \|\psi\|_{\infty}\le
|A|= A_b(\varphi)$, and we have finished. {\hskip 5pt}
$\blacksquare$
\renewcommand{\qedsymbol}{}
\end{proof}
We also have the following result for ${\rm BMO}_{2s}$.
\begin{proposition} A function $\varphi\in {\rm BMO}$ if and only if
$\varphi\in {\rm BMO}_{2s}$ and
\[A_0(\varphi)=\sup_n \Big(\,\frac{1}{2^{n+1}}\Big|\int_{[-2^n,0]}
\varphi(y)\,dy - \int_{[0,2^n]} \varphi(y)\,dy\Big|\,\Big)
<\infty\,.\] Moreover, there is a constant $c$ such that
\[\|\varphi\|_*\sim \max(\,\|\varphi\|_{*,2s},A_0(\varphi)\,).\]
\end{proposition}
We leave the verification of this fact to the reader, and point
out an interesting consequence, see \cite{Gerard}.
\begin{proposition}Suppose $\varphi\in {\rm BMO}_{2s}$ is  supported
in $[\,0,\infty)$.
\begin{enumerate}[\upshape 1.]
\item The even extension $\varphi_e$ of $\varphi$ belongs to {\rm
BMO}, and $\|\varphi_e\|_*=\|\varphi\|_{*,2s}$. \item  The odd
extension $\varphi_o$ of $\varphi$ belongs to {\rm BMO} if and
only if
\[\sup_n\frac{1}{2^n}\Big|\int_{[\,0,2^n]}\varphi(y)\,dy\Big|<\infty\,,\]
and in this case \[\|\varphi_o\|_*\sim \|\varphi\|_{*,2s}+
\sup_n\frac{1}{2^n}\Big|\int_{[\,0,2^n]}\varphi(y)\,dy\Big|.\]
\end{enumerate}
\end{proposition}
 As an illustration of the use of the above results we
 will consider the $T(1)$ Theorem, which
 establishes the continuity in $L^2$ of a standard CZO operator
 essentially  under two kinds of  assumptions, the weak boundedness
property  and the $T(1), T^*(1)$ BMO assumption. Indeed, we have,
see \cite{davidjourne},
\begin{theorema}
Suppose $T$ is a standard {\rm CZO} that satisfies
\begin{enumerate}[]
\item {\rm (WBP) } For every interval $I$,\, $|\langle
T\chi_I,\chi_I\rangle|\le \, c\,|I|\,.$ \item {\rm (BMO}
condition{\rm ) } $\|T(1)\|_*+ \|T^*(1)\|_*\le\,c\,.$
\end{enumerate}
Then $T$ is a continuous mapping in $L^2$.
\end{theorema}
In applications it is of interest to state the BMO condition in a
form that is easily verified. For instance, in the dyadic setting,
the following two conditions may be assumed instead,
\begin{enumerate} \item[${\rm (1}_d{\rm )}$] $\|T(1)\|_{*,d}+
\|T^*(1)\|_{*,d}\le\,c\,.$ \item[${\rm (2}_d{\rm )}$] For all
integers $n,k$, $ |\langle T(b_{n,k}),1\rangle| +\, |\langle
T^*(b_{n,k}),1\rangle| \le \,c\,.$
\end{enumerate}
Then clearly $T(1)$, $T^*(1)\in {\rm{BMO}}$, and the $T(1)$
Theorem obtains.

Similarly, in the two-sided setting, the following two conditions
may be used instead of the BMO assumption,
\begin{enumerate} \item[${\rm (1}_s{\rm )}$]
$\|T(1)\|_{*,2s}+ \|T^*(1)\|_{*,2s}\le\,c\,.$ \item[${\rm
(2}_s{\rm )}$] For all integers $n$, $ |\langle
T(b_{n,0}),1\rangle| +\, |\langle T^*(b_{n,0}),1\rangle| \le
\,c\,.$
\end{enumerate}

 Two  particular instances of this last observation
  come to mind. Let $T$ be a CZO with WBP
 that  satisfies ${\rm (1}_s{\rm )}$. Also, assume that for any interval
$I=I_{n,0}$,  $T(\chi_I)$ is supported in $I$, and similarly for
$T^*$. Now, since
\[\langle T(b_{n,0}),1\rangle=\frac1{2^{n+1}}\left[\, \langle\,
T\chi_{[-2^n,0]}, \chi_{[-2^n,0]} \,\rangle -\langle\,
T\chi_{[0,2^n]}, \chi_{[0,2^n]}\,\rangle  \,\right]\!,
\]
by WBP, $|\langle T(b_{n,0}),1\rangle|\le c$. The estimate for
$T^*$ is obtained in a similar fashion, and  therefore ${\rm
(2}_s{\rm )}$ also holds. Thus $T$ is bounded in $L^2$.

Finally, when the kernel of $T$ is even, or odd, in $x$ and $y$,
 $T(1)$ and $T^*(1)$ are  even, or odd,
respectively. Now, by Proposition 3.4, if $T(1)$ is even and
$T(1)\,\chi_{[0,\infty)}$ is in ${\rm BMO}_{2s}$, then $T(1)\in
{\rm BMO}$; similarly for $T^*(1)$. On the other hand, if $T(1),
T^*(1)$ are odd and $T(1)\,\chi_{[0,\infty)},
T^*(1)\,\chi_{[0,\infty)}$ are in ${\rm BMO}_{2s}$,  we also
require that
\[\sup_n\frac{1}{2^n}\Big|\int_{[0,2^n]}T(1)(y)\,dy\,\Big|,
\ \sup_n\frac{1}{2^n}\Big|\int_{[0,2^n]}T^*(1)(y)\,dy\,\Big| \le
c\,.\] Under these assumptions ${\rm (1}_s{\rm )}$ holds and
together with ${\rm (2}_s{\rm )}$ obtain the continuity of $T$ in
$L^2$.
\section{Final remarks}
We  sketch  now  the extension of the results to higher
dimensions. To avoid  technicalities  we restrict ourselves to the
case $n=2$, but stress that appropriate versions remain valid for
 arbitrary $n$. Also, since the proofs  follow along similar lines
 to the case $n=1$, they will be omitted.

 The two-dimensional Haar system is generated by the integer
 translations and  dyadic dilations  of the three basic
orthogonal functions
\[
\Psi_1(x,y) = H(x) \, \chi_{[0,1]} (y)\,,\quad \Psi_2(x,y) =
\chi_{[0,1]} (x) \, H(y)\,,\]  \[\Psi_3(x,y) = H(x) \, H(y)\,,
\]
cf. \cite{Wojtaszszyk}.  More precisely, the functions
\[\Psi_{1,n,k,l}(x,y)=2^n \Psi_1(2^n x- k,2^ny-l)\,,
\Psi_{2,n,k,l}(x,y)=2^n\Psi_2(2^nx-k,2^ny-l) \,,\] \[
\Psi_{3,n,k,l}(x,y)=2^n\Psi_3(2^nx-k,2^ny-l)\,,\] for arbitrary
integers $n,k,l$, generate the two dimensional Haar system. In
three dimensions seven basic functions are required.

We then have
\begin{theorem} The closed span of the two dimensional
Haar system in  $H^1(R^2)$  is the subspace $H^A(R^2)$ of
$H^1(R^2)$ which consists of those functions that have $0$
integral on each quadrant.
\end{theorem}

Clearly there is some redundancy in this statement. Since
functions in $H^1(R^2)$ have $0$ integral, it suffices to require
that the functions in question have $0$ integral in any three
quadrants.

Let $Q_1$, $Q_2$, $Q_3$, $Q_4$, denote the four quadrants of
$R^2$. It is not hard to see that the functions $\varphi$ with the
property that $\int\int f\varphi=0$ for all $f\in{\rm
sp}\{\Psi_{j,n,k}\}$ are of the form $\varphi=\sum_i
c_i\,\chi_{Q_i}$ and this suggests how close $H^A(R^2)$ is to
$H^1(R^2)$. In fact, we have,
\[H^1(R^2)=\overline{\rm {sp}}\,\{\Psi_{j,n,k},b_1,b_2,b_3\}\,,\] where
 $b_1(x,y)=b(x)\chi_{[0,1)}(y)$, $b_2(x,y)=\chi_{[0,1]}(x)b(y)$,
 and $b_3(x,y)=$\\ $b(x)\,\chi_{[-1,0]}(y)$.

As for arbitrary atoms $a$ in $H^1(R^2)$,  they   can be expressed
as a sum of at most five atoms, four dyadic and a special atom.
More precisely, if we denote $Q_{n,k,m,l}=I_{n,k}\times I_{m,l}$,
then
\begin{lemma}
 Let $a$ be an $H^1(R^2)$ atom. Then there are at most
five atoms $a_1, a_2, a_3, a_4, b_{n,k,m,l}$, such that
\begin{enumerate}[\upshape i.]
\item The $a_i\!$'s are dyadic atoms. \item For some integers
$n,k,m,l$, $b_{n,k,m,l}(x)$ is equal to
\[ \frac {1}{2^{n+m+2}}\left[ k_1\chi_{Q_{n,k,m,l}}
+k_2\chi_{Q_{n,k-1,m,l}} +k_3 \chi_{Q_{n,k,m,l-1}}+ k_4
\chi_{Q_{n,k-1,m,l-1}}\right],\] where $|k_i|\le c$, an absolute
constant independent of $a$, and $\sum_1^4 k_i=0$. \item
$a=\sum_1^4c_ia_i+c_5b_{n,k,m,l}$, where $|c_i|\le c$, an absolute
constant independent of $a$.
\end{enumerate}
\end{lemma}
With this decomposition of individual atoms available, we can
describe $H^1(R^2)$ in various ways. For instance, if $H^1_d(R^2)$
denotes the dyadic Hardy space, and $H^1_{\delta}(R^2)$ denotes
the subspace of $H^1(R^2)$ spanned by the $b_{n,k,m,l}$'s, then
$H^1(R^2)=H^1_d(R^2) + H^1_{\delta}(R^2)$ in the sense of sum of
Banach spaces. Thus the continuity of a linear operator acting on
$H^1(R^2)$ can be characterized in terms of the continuity of its
restrictions to $H^1_d(R^2)$ and $H^1_{\delta}(R^2)$.

Also, decompositions of $H^1(R^2)$ lead to characterizations of
BMO$(R^2)$. More precisely,
\begin{lemma} A locally integrable $\varphi\in{\rm BMO}(R^2)$ if
and only if $\varphi$ belongs to ${\rm BMO}_d(R^2)$ and for a
constant $c$,
\[\Big|\int\int_{R^2}\varphi(x,y)\,b_{n,k,m,l}(x,y)\,dx\,dy\,\Big|
\le \, c\,,\quad  {\rm{\,all}}\, n,k,m,l\,.\]
\end{lemma}

Finally, as a consequence of this result one can write a version
of the $T(1)$ theorem in $R^2$ under dyadic-like assumptions.

DEPARTMENT OF MATHEMATICS, INDIANA UNIVERSITY, BLOOMINGTON, IN
47405
\\{\it E-mail:} wabusham@indiana.edu, torchins@indiana.edu

\begin{thebibliography}{99}

\bibitem{Gerard} {G. Bourdaud},
{\em Remarques sur certains sous-espaces de ${\rm BMO}\,(R^n)$ et
de ${\rm bmo}\,(R^n)$},  Ann. Inst. Fourier (Grenoble) {\bf 52}
(2002), no. 4, 1187-1218.

\bibitem {bownik}  {M. Bownik}, {\em Boundedness of operators on Hardy
spaces via atomic decompositions},  Proc. Amer. Math. Soc. {\bf
133} (2005), no. 12, 3535--3542.

\bibitem{APC} {A. P. Calder\'{o}n}, {\em Estimates for singular
integral operators in terms of maximal functions}, Studia Math
{\bf 44} (1972), 563--582.

\bibitem{APCAT} {A. P. Calder\'{o}n} and {A. Torchinsky},
{\em Parabolic maximal functions associated with a distibution,
II}, Advances in Math., {\bf 249} (1977), 101--171.

\bibitem{Coifman} {R. R. Coifman},
{\em A real variable characterization of $H^p$},  Studia Math.
{\bf 51} (1974), 269--274.

\bibitem{davidjourne} {G. David} and  {J. Journ\'{e}}, {\em A boundedness
criterion for generalized Calder\'{o}n-Zygmund operators}, Ann. of
Math. {\bf 120} (1984), 371--397.

\bibitem{desouza} {G. S. de Souza}, {\em Spaces formed by special atoms, I, }
 Rocky Mountain J.  Math. {\bf 14} (1984),  no.2, 423--431.

\bibitem{Fefferman_Stein} {C. Fefferman} and {E. M. Stein},
{\em $H^p$ spaces of several variables},  Acta. Math. {\bf 129}
(1972), 137--193.

\bibitem{fridli} {S. Fridli}, {\em Transition from the dyadic to the real nonperiodic
Hardy space},  Acta Math. Acad. Paedagog. Nih\'{a}zi (N.S.) {\bf
16} (2000), 1--8, (electronic).

\bibitem{fridli1}  {S. Fridli}, {\em Hardy spaces generated by an integrability
condition},  J. Approx. Theory, {\bf 113} (2001), no.1, 91--109.

\bibitem{GCRdF} {J. Garc\'{i}a-Cuerva} and {J. L. Rubio de
Francia}, {\em Weighted norm inequalities and related topics},
Notas de Matem\'{a}tica {\bf 116}, North Holland, 1985.

\bibitem{garnettjones} {J. Garnett} and {P. Jones}, {\em ${\rm BMO}$ from dyadic
${\rm BMO}$},  Pacific J. Math. {\bf 99} (1982),  no. 2, 351--371.

\bibitem{Haar} {A. Haar},
{\em Zur Theorie der orthogonalen Funktionensysteme},  Math. Ann.
{\bf 69} (1910), 331--371.

\bibitem{JN} {F. John} and {L. Nirenberg}, {\em On functions of bounded mean
oscillation}, Comm. Pure Appl. Math. {\bf 14} (1961), 415--426.

\bibitem{Latter} {R. H. Latter}, {\em A characterzation
of $H^p(R)$ in terms of atoms},
 Studia Math.  {\bf 62} (1978), 93--101.

\bibitem {tao_mei} {T. Mei}, {\em {\rm BMO} is the intersection of two
translates of dyadic {\rm BMO}}, C. R. Math. Acad. Sci. Paris {\bf
336} (2003),  no. 12, 1003--1006.

\bibitem{Myer_English} {Y. Meyer}, Wavelets and operators,
{\em Cambridge University Press, Cambridge}, 1992.

\bibitem{shiu} {J.-L. Shiu}, {\em The $H^1$-closure of the Haar system
and its dual space},  Ph D dissertation, Indiana University, 2004.

\bibitem{stein} {E. M. Stein}, {\em Harmonic analysis: Real-variable methods,
orthogonality and oscillatory integrals},  Princeton Math Ser.
{\bf 43}, Princeton  University Press, Princeton, 1993.

\bibitem{torchinsky} {A. Torchinsky},
{\em Real-variable methods in harmonic analysis},  Dover
Publications, Inc., 2004.

\bibitem{Wojtaszszyk} {P. Wojtaszczyk}, {\em A mathematical
introduction to
 wavelets}, London Mathematical Society Student Text, {\bf 37},  Cambridge
 University Press, Cambridge, 1997.
\end{thebibliography}
\end{document}